\documentclass[11pt]{amsart}
\usepackage{amsthm,amsbsy,amsfonts,amssymb,amsmath,amscd}

\newcommand{\Q}{\mathbb Q}
\newcommand{\Z}{\mathbb Z}
\newcommand{\R}{\mathbb R}
\newcommand{\C}{\mathbb C}

\newcommand{\hH}{\mathcal H}
\newcommand{\nN}{\mathcal N}
\newcommand{\A}{\mathbb A}
\newcommand{\I}{\mathcal A}

\newcommand{\s}{\text sym}
\makeatletter
\def\l@section{\@tocline{1}{4pt}{1pc}{}{}}
\def\l@subsection{\@tocline{2}{0pt}{2pc}{5pc}{}}
\makeatother

\title{Remarks on the symmetric powers of  cusp forms on GL$(2)$}

\author{Dinakar Ramakrishnan}

\begin{document}
\maketitle

\medskip

\begin{flushright}
{\it To Steve Gelbart \\ On the occasion of his sixtieth birthday}
\end{flushright}

\bigskip

\section*
{\bf Introduction}

\bigskip

Let $F$ be a number field, and $\pi$ a cuspidal automorphic
representation of GL$(2,\A_F)$ of conductor $\nN$. For every $m\geq
1$ one has its {\it symmetric $m$-th power $L$-function}
$L(s,\pi;{\rm sym}^m)$, which is an Euler product over the places
$v$ of $F$, with the $v$-factors (for finite $v\nmid \nN$ of norm
$q_v$) being given by
$$
L_v(s,\pi;{\rm sym}^m) \, = \, \prod_{j=0}^m
(1-\alpha_v^j\beta_v^{m-j}{q_v}^{-s})^{-1},
$$
where the unordered pair $\{\alpha_v, \beta_v\}$ defines the
diagonal conjugacy class in GL$_2(\C)$ attached to $\pi_v$. Even at
a ramified (resp. archimedean) place $v$, one has by the local
Langlands correspondence a $2$-dimensional representation $\sigma_v$
of the extended Weil group $W_{F_v}\times{\rm SL}(2,\C)$ (resp. of
the Weil group $W_{F_v}$), and the $v$-factor of the symmetric
$m$-th power $L$-function is associated to sym$^m(\sigma_v)$. A special case of the
{\it principle of functoriality} of Langlands asserts that there is, for each $m$,
an (isobaric) automorphic representation ${\rm sym}^m(\pi)$ of
GL$(m+1,\A)$ whose standard (degree $m+1$) $L$-function $L(s, {\rm
sym}^m(\pi))$ agrees, at least at the primes not dividing $\nN$,
with $L(s,\pi;{\rm sym}^m)$. It is well known that such a result
will have very strong consequences, such as the {\it Ramanujan
conjecture} and the {\it Sato-Tate conjecture} for $\pi$. The {\it
modularity}, also called {\it automorphy}, has long been known for
$m=2$ by the pioneering work of Gelbart and Jacquet (\cite{GeJ}); we
will write Ad$(\pi)$ for the selfdual representation
sym$^2(\pi)\otimes\omega^{-1}$, $\omega$ being the central character
of $\pi$. A {\it major breakthrough}, {\sl due to Kim and Shahidi}
(\cite{KSh1, KSh2, K})), has established the modularity of
sym$^m(\pi)$ for $m=3,4$, along with a useful cuspidality criterion
(for $m\leq 4$). Furthermore, when $F=\Q$ and $\pi$ is defined by a
holomorphic newform $f$ of weight $2$, $\Q$-coefficients and level
$N$, such that at some prime $p$, the component $\pi_p$ is
Steinberg, a recent {\it dramatic theorem} of Taylor, Harris, Clozel 
and Shepherd-Baron, furnishes the {\it potential modularity}
of sym$^{2m}(\pi)$ (for every $m\geq 1$), i.e., its modularity over
a number field $K$, thereby extracting the Sato-Tate conjecture in
this case by a clever finesse. It should however be noted that such a beautiful result
is not (yet) available for $\pi$ defined by newforms $\varphi$ of
higher weight, for instance for the ubiquitous cusp form $\Delta(z)
\, = \, q\prod_{n \geq 1} (1-q^n)^{24} \, = \, \sum_{n \geq 1}
\tau(n) q^n$, where $z\in\hH$ and $q = e^{2\pi iz}$, which is
holomorphic of weight $12$, level $1$ and trivial character.

\medskip

In this Note we consider the following more modest, but nevertheless
basic, question:

\medskip

{\it Suppose sym$^m(\pi)$ is an automorphic representation of
GL$_{m+1}(\A_F)$. When is it cuspidal?}

\medskip

If sym$^m(\pi_v)$ is, for some finite place $v$, in the discrete
series, which happens for example when $\pi_v$ is Steinberg, it is
well known that the global representation sym$^m(\pi)$ will
necessarily be cuspidal (once it is automorphic). On the other hand,
one knows already
for $m=2$, as shown by Gelbart and Jacquet (\cite{GeJ}), that if $\pi$ is
dihedral, i.e., associated to an idele class character $\chi$ of a
quadratic extension $K$ of $F$, then sym$^2(\pi)$ is not cuspidal; in
fact, this is {\it necessary and
sufficient} condition. There is a
non-trivial extension of such a criterion in the work of Kim and
Shahidi (\cite{KSh2}), who show that for a non-dihedral $\pi$,
sym$^3(\pi)$ is Eisensteinian iff $\pi$ is {\it tetrahedral}, while
sym$^4(\pi)$ is cuspidal iff $\pi$ is not tetrahedral or octahedral.
We will say that $\pi$ is {\it solvable polyhedral} iff it is
dihedral, tetrahedral or octahedral. Finally, if $\pi$ is associated
to an irreducible $2$-dimensional Galois representation $\rho$ which
is icosahedral, i.e., with projective image isomorphic to the
alternating group $A_5$, one knows that sym$^6(\rho)$ is reducible,
suggesting that sym$^6(\pi)$ is not cuspidal. However, sym$^5(\rho)$ is, in the icosahedral case,
necessarily a tensor product sym$^2(\rho')\otimes\rho$, where
$\rho'$ is the Galois conjugate representation of $\rho$ (which is
defined over $\Q[\sqrt{5}]$) (cf. \cite{K2}, \cite{Wang}, for example). This allowed Wang to prove
(in \cite{Wang}) that
sym$^5(\pi)$ is cuspidal by making use of the construction (cf
\cite{KSh1}) of the functorial product $\Pi\boxtimes\pi^\prime$ (in
GL$(6)/F$), for $\Pi$ (resp. $\pi^\prime$) a cusp form on GL$(3)/F$
(resp. GL$(2)/F$), and by developing a cuspidality criterion for this product.

\medskip

In order to answer the question, we make the following definition:
Call an irreducible cuspidal automorphic representation $\pi$ of GL$(2,\A_F)$
{\it quasi-icosahedral} iff we have
\begin{enumerate}
\item[(i)]sym$^m(\pi)$ is automorphic for every $m \leq 6$;
\item[(ii)]sym$^m(\pi)$ is cuspidal for every $m \leq 4$;
and
\item[(iii)]sym$^6(\pi)$ is not cuspidal.
\end{enumerate}
The key result which we prove (see part (b) of Theorem A below) is that, for every such quasi-icosahedral $\pi$ of central character $\omega$,
there exists another
cusp form $\pi'$ of GL$(2)/F$ (of central character $\omega'$) such that the symmetric fifth power
of such a quasi-icosahedral cusp form $\pi$ is necessarily a character twist of the functorial
product ${\rm Ad}(\pi')\boxtimes\pi$, where Ad$(\pi')=$sym$^2(\pi')\otimes{\omega'}^{-1}$. If $\pi$ were associated to an icosahedral Galois representation $\rho$,
defined over $\Q[\sqrt{5}]$, $\pi'$ could be taken to correspond to the Galois conjugate
representation $\rho^{[\theta]}: \sigma\to\rho(\theta\sigma\theta^{-1})$, where $\theta$ denotes the non-trivial
automorphism of $\Q[\sqrt{5}]$. The beauty is that we can find $\pi'$ by a purely automorphic argument.

\medskip

All of this is consistent with the results of Wang , as well as
with the philosophy of Langlands (\cite{La}), which
predicts that any cuspidal $\pi$ on GL$(2)/F$
should be naturally associated to a reductive subgroup $H(\pi)$ of
GL$_2(\C)$. In fact, one expects there to be a pro-reductive group ${\mathcal L}_F$
over $\C$ whose $n$-dimensional $\C$-representations $\sigma$ classify (up to equivalence)
the (isobaric) automorphic representations $\pi$ of GL$(n,\A_F)$,
and $H(\pi)$ should be given by the image of $\sigma$.

\medskip

\noindent{\bf Theorem A} \, \it Let $\pi$ a cuspidal automorphic
representation of GL$_2(\A_F)$, which is not solvable polyhedral, of
central character $\omega$. Suppose sym$^m(\pi)$ is modular for all
$m$. Then we have
\begin{enumerate}
\item[(a)] sym$^5(\pi)$ is cuspidal.
\item[(b)] sym$^6(\pi)$ is non-cuspidal iff we have
$$
{\rm sym}^5(\pi) \, \simeq \, {\rm Ad}(\pi')\boxtimes
\pi\otimes\omega^2,
$$
for a cuspidal automorphic representation $\pi'$ of GL$_2(\A_F)$.
\item[(c)] If sym$^6(\pi)$ is cuspidal, then so is sym$^m(\pi)$ for
all $m\geq 1$.
\item[(d)] If $F=\Q$ and $\pi$ is
defined by a non-CM, holomorphic newform $\varphi$ of weight $k\geq
2$, then ${\rm sym}^m(\pi)$ is cuspidal for all $m$.
\end{enumerate}
\rm

\medskip

One can do a bit better than this in that for a given symmetric
power, one does not need information on all the sym$^m(\pi)$. See
Theorem A$'$ in section 2 for a precise statement. The proofs are then
given in sections 3 and 4.

\medskip

In part (b), the cusp form $\pi'$ is not uniquely determined, only up to
a character twist. In a sequel we will show that, in fact, for a suitable
choice of $\pi'$, ${\rm sym}^5(\pi)$ is, in the quasi-icosahedral case,
expressible as a character twist of
${\rm Ad}(\pi)\boxtimes \pi'$; $\pi'$ will also turn out to be quasi-icosahedral.
This is as predicted by looking at the
Galois side, and it will help us normalize the choice of $\pi'$,
leading in addition to a precise rationality statement.

\medskip

The results of this paper were essentially established some time
ago, but the questions raised to me in the past two years by some
colleagues have led me to believe in the possible usefulness of
their being in print. While the inspiration for the results here
came from Langlands (and the paper of Wang), and from a short conversation
with Richard Taylor some time back, the proofs depend, at
least partly, on the beautiful constructions \cite{KSh1, KSh2, K} of
Kim and Shahidi. Use is also made of the papers \cite{Ra2, Ra6}.

\medskip

\noindent{\bf Acknowledgement}: Like so many others interested in
Automorphic Forms, I was decidedly influenced during my graduate
student years (in the late seventies), and later, by Steve Gelbart's
book, {\it Automorphic Forms on adele groups}, and his expository
papers, {\it Automorphic forms and Artin's conjecture} and {\it
Elliptic curves and automorphic representations}, as well as his
seminal work with Jacquet, {\it A relation between automorphic forms
on ${\rm GL}(2)$ and ${\rm GL}(3)$}. His later works have also been
influential. Furthermore, Steve has been very friendly and
generous over the years, and it is a great pleasure to dedicate this
paper to him. I would like to thank Freydoon Shahidi and Erez Lapid
for their interest in this paper, and especially
the latter for reading the manuscript in detail
and making comments, resulting in a streamlining of the exposition.
Finally, I would be remiss if I do not acknowledge
support from the NSF through the grant DMS0701089.

\vskip 0.2in

\section{Preliminaries}

\bigskip

\subsection{The standard $L$-function of GL$(n)$}

Let $F$ be a number field with adele ring $\A_F$. For each place $v$, denote by
$F_v$ the corresponding local completion of $F$,
and for $v$ finite, by $\mathfrak O_v$ the ring of integers of $F_v$ with uniformizer
$\varpi_v$ of norm $q_v$.
For any algebraic group $G$ over $F,$ let $G(\A_F)$ denote the
restricted direct product $\prod'_v G(F_v),$ endowed with the
usual locally compact topology.
For $m\geq 1$, let $Z_m$ denote the center of GL$(m)$.
One knows that the volume of
$Z_m(\A_F){\rm GL}_m(F)\backslash {\rm GL}(m, \A_F))$ is finite.

By a {\it unitary cuspidal} representation of GL$_m(\A_F) =
{\rm GL}_m(F_\infty) \times {\rm GL}_m(\A_{F,f})$, we will always mean an
irreducible, automorphic representation occurring in the space of
cusp forms in $L^2(Z_m(\A_F)G_m(F)\backslash G_m(\A_F), \omega)$
relative to a character $\omega$ of $Z_m(\A_F)$,
trivial on $Z_m(F)$. By a (general) {\it
cuspidal representation} of GL$_m(\A_F)$, we will mean an irreducible admissible
representation of $G_m(\A_F)$ for which there exists a real
number, called the {\it weight of $\pi$} such that $\pi \otimes
|.|^{w/2}$ is a unitary cuspidal representation. Such a
representation is in particular a restricted tensor product $\pi
\, = \, \otimes'_v \pi_v \, = \, \pi_{\infty} \otimes \pi_f$,
where each $\pi_v$ is an (irreducible) admissible representation
of $G(F_v)$, with $\pi_v$
unramified at almost all $v.$

For any irreducible, automorphic representation $\pi$ of
$GL(n,\A_F),$ let $L(s, \pi)$ $= L(s, \pi_{\infty})L(s, \pi_f)$ denote
the associated {\it standard} $L-$function (\cite{J}) of $\pi;$ it
has an Euler product expansion
$$
L(s,\pi) \, = \, \prod_v \, L(s, \pi_v), \leqno (1.1.1)
$$
convergent in a right-half plane. If $v$ is an archimedean place,
then one knows (cf. \cite{La1}) how to associate a semisimple
$n-$dimensional $\C-$representation $\sigma(\pi_v)$ of the Weil
group $W_{F_v},$ and $L(\pi_v,s)$ identifies with $L(\sigma_v,s).$
On the other hand, if $v$ is a finite place where $\pi_v$ is
unramified, there is a corresponding semisimple (Langlands)
conjugacy class $A_v(\pi)$ (or $A(\pi_v)$) in GL$(n,\C)$ such that
$$
L(s,\pi_v) \, = \, {\rm {det}}(1-A_v(\pi)T)^{-1}|_{T=q_v^{-s}}.
\leqno (1.1.2)
$$
We may find a diagonal representative diag$(\alpha_{1,v}(\pi), ...
, \alpha_{n,v}(\pi))$, unique up to
permutation of the diagonal entries, for $A_v(\pi)$ . Let $[\alpha_{1,v}(\pi), ...
, \alpha_{n,v}(\pi) ]$ denote the resulting unordered $n-$tuple.
Since $W_{F,v}^{{\rm {ab}}} \simeq F_v^\ast,$ $A_v(\pi)$ clearly
defines an abelian $n-$dimensional representation $\sigma(\pi_v)$
of $W_{F,v}.$ One has

\medskip

\noindent{\bf Theorem 1.1.3} (\cite{GoJ, J}) \quad \it Let $n \geq
1,$ and $\pi$ a cuspidal representation of GL$(n,\A_F)$. When $n=1$,
assume that $\pi$ is not of the form $\vert\cdot\vert^w$ for any $w\in \C$.
Then $L(s,\pi)$ is entire. \rm

\medskip

When $n = 1,$ such a $\pi$ is simply a unitary idele class
character $\chi$, and this result is due to Hecke. Also, when $\chi$ is
trivial, $L(s,\pi_f)$ is the Dedekind zeta function $\zeta_F(s).$

\medskip

Given a pair of automorphic representations $\pi, \pi'$ of GL$(n,\A_F)$, GL$(n',\A_F)$,
respectively, one can associate an $L$-function $L(s, \pi\times \pi')$ which is meromorphic.
We will postpone its definition till section 1.4.

\medskip

For any $L$-function with an Euler product expansion (over $F$): $L(s)=\prod_v L_v(s)$,
if $S$ is any set of places of $F$, the associated
{\it incomplete $L$-function} is defined as follows:
$$
L^S(s) : = \, \prod_{v\notin S} L_v(s).
$$

\bigskip

\subsection{Isobaric automorphic representations}

By the theory of Eisenstein series, one has a sum operation
$\boxplus$ (\cite{La3}) on a suitable set of automorphic representations
of GL$(n)$ for all $n$. One has the following:

\medskip

\noindent{\bf Theorem 1.2.1} \, ([JS])\quad \it Given any
$m-$tuple of cuspidal representations $\pi_1, ..., \pi_m$ of
GL$(n_1,\A_F), ... ,$ GL$(n_m,\A_F)$ respectively, there exists an
irreducible, automorphic representation $\pi_1 \boxplus ...
\boxplus \pi_m$ of GL$(n,\A_F),$ $n \, = \, n_1 + ... + n_m$, which is unique and satisfies the
following property: For any finite set $S$ of
places, we have, for every cuspidal automorphic representation $\pi'$ of GL$(n',\A_F)$
(with $n'$ arbitrary),
$$
L^S(s, \left(\boxplus_{j=1}^m \pi_j\right) \times \pi') \, = \, \prod_{j=1}^m L^S(s,
\pi_j\times \pi').
$$
\rm

\medskip

The $L$-functions in the Theorem are the Rankin-Selberg $L$-functions attached to pairs
of automorphic representations, which we briefly discuss in section 1.4 below.

\medskip

Call such a (Langlands) sum $\pi \simeq \boxplus_{j=1}^m \pi_j$,
with each $\pi_j$ cuspidal, an {\it isobaric} representation. Denote
by ram$(\pi)$ the finite set of finite places where $\pi$ is
ramified, and let $\mathfrak N(\pi)$ be its conductor
(\cite{JPSS2}).

For every integer $n \geq 1,$ set:
$$
\mathcal  A(n,F) \, = \, \{\pi: {\rm {isobaric}} \, {\rm {
representation}} \, {\rm {of}} \, {\rm {GL}}(n,\A_F) \}/{\simeq },
\leqno (1.2.2)
$$
\noindent and
$$
\mathcal  A_0(n,F) \, = \, \{\pi \in \mathcal  A(n,F) | \, \pi \,
{\rm {cuspidal}} \}.
$$
Put $\mathcal  A(F) \, = \, \cup_{n \geq 1} \mathcal  A(n,F)$ and
$\mathcal  A_0(F) \, = \, \cup_{n \geq 1} \mathcal  A_0(n,F).$

\medskip

\noindent{\bf Definition 1.2.3} \, \it Given $\pi, \eta\in{\mathcal
A}(F)$, if we can find an $\eta' \in {\mathcal A}(F)$ such that $\pi
\, \simeq \, \eta \boxplus \eta'$, we will call $\eta$ an isobaric
summand of $\pi$ and write
$$
[\pi : \eta] \, > \, 0.
$$
\rm

\medskip

\noindent{\bf Remark}. \quad One can also define the analogues of
$\mathcal A(n, F)$ for local fields $F$, where the ``cuspidal''
subset $\mathcal A_0(n, F)$ consists of essentially
square-integrable representations of GL$(n, F)$. See \cite{La3} and
\cite{Ra1} for details.

\bigskip

\subsection{Symmetric powers of GL$(2)$}

\medskip

Since the $L$-group of GL$(2)$ is GL$(2, \C) \times W_F$, the
principle of functoriality of Langlands (\cite{La2}) predicts that
for any algebraic representation
$$
r: {\rm GL}(2, \C) \, \rightarrow \, GL(N, \C), \leqno(1.3.1)
$$
and any number field $F$, there should be a map
$$
{\mathcal A}(2,F) \, \rightarrow \, {\mathcal A}(N,F), \, \pi \to
r(\pi), \leqno(1.3.2)
$$
with compatible local maps, such that for all finite unramified
places $v$ (for $\pi$), we have the equality of Langlands classes
$$
r(A(\pi_v)) \, = \, A(r(\pi)_v).
$$

It suffices to establish this for irreducible representations $r$,
which are all of the form ${\rm sym}^n(r_0) \otimes L^{\otimes k}$,
with $n, k \in \Z, n\geq 0$; here $r_0$ denotes the
standard representation of GL$(2, \C)$ with determinant $L$, and
sym$^n(r_0)$ denotes the symmetric $n$-th power representation of
$\rho$.

It is enough to construct the sym$^n(\pi)$'s for $\pi$ cuspidal. When
it exists, by which we mean it belongs to ${\mathcal A}(F)$, we will
write (for $\pi \in {\mathcal A}(2,F)$)
$$
{\rm sym}^n(\pi) \, = \, {\rm sym}^n(\rho)(\pi).
$$
It may be useful to recall that if
$$
L(s, \pi_v) \, = \, [(1-\alpha_vq_v^{-s})(1-\beta_vq_v^{-s})]^{-1}
$$
at any unramified finite place $v$ with norm $q_v$, with $A(\pi_v)$
being represented by the diagonal matrix with entries $\alpha_v,
\beta_v$, then for every $n \geq 1$,
$$
L(s, \pi_v, {\rm sym}^n) \, = \, [\prod_{j=0}^{n}
(1-\alpha_v^j\beta_v^{n-j}q_v^{-s})]^{-1}.\leqno(1.3.3)
$$

It is well known that when $r = L$, $r(\pi) \in {\mathcal A}(1, F)$
is given by the central character $\omega = \omega_\pi$ of $\pi$.
Consequently, if one can establish
the lifting for $r = {\rm sym}^n(r_0)$, then one can also achieve
it for $r = {\rm sym}^n(r_0) \otimes L^{\otimes k}$ by twisting by
$\omega^k$, i.e., by setting
$$
\left({\rm sym}^n(r_0) \otimes L^{\otimes k}\right )(\pi) \, = \,
{\rm sym}^n(\pi) \otimes \omega^k.
$$
So it suffices to establish the transfer $\pi \to r(\pi)$ for
sym$^n(r_0)$ for all $n$. Clearly, ${\rm sym}^1(\pi) \, = r_0(\pi) \, = \, \pi$.

\medskip

\noindent{\bf Proposition 1.3.4} \, \it Let $\pi$ be a cuspidal
automorphic representation of GL$(2, \A_F)$ which is associated to a
two-dimensional, continuous $\C$-representation $\rho$ of
Gal$(\overline F/F)$ so that $L(s, \rho) = L(s, \pi)$. Suppose ${\rm
sym}^m(\pi)$ exists in ${\mathcal A}(F)$ for every $m \geq 1$. It is
then cuspidal iff ${\rm sym}^m(\rho)$ is irreducible. \rm

\medskip

\noindent{\it Proof}. \,
For any continuous finite-dimensional $\C$-representation $\sigma$ of
$\Gamma_F:=$Gal$(\overline F/F)$, one knows (cf. \cite{Tate}) the following fact about
the Artin $L$-functions:
$$
{\rm Hom}_{\Gamma_F}(\C, \sigma) \, = \,
-{\rm ord}_{s=1}L^S(s, \sigma).\leqno(a)
$$
Applying this to
$$
\sigma:={\rm sym}^m(\rho)\otimes{\rm sym}^m(\overline \rho) \, \simeq \, {\rm End}(\rho),
$$
we see by Schur's lemma that
$$
{\rm sym}^m(\rho) \, \, irreducible \quad \Longleftrightarrow
-{\rm ord}_{s=1}L^S(s, {\rm sym}^m(\rho)\otimes {\rm sym}^m({\overline \rho})) = 1.\leqno(b)
$$
On the other hand, by a result we will prove later in section 3.1 (see Lemma 3.1.1), sym$^m(\pi)$ is,
when it is automorphic,
an isobaric sum of {\it unitary} cuspidal representations. This implies that
$$
{\rm sym}^m(\pi) \, \, cuspidal \quad \Longleftrightarrow
-{\rm ord}_{s=1}L^S(s, {\rm sym}^m(\pi)\times {\rm sym}^m({\overline \pi})) = 1,\leqno(c)
$$
where the $L$-function is the Rankin-Selberg $L$-function (see 1.3 below for its basic
properties).
Finally, since by hypothesis, $\rho$ corresponds to $\pi$, the $L$-functions of (b) and
(c) are the same. The assertion follows.

\qed

\medskip

One expects the same when $\rho$ is an $\ell$-adic Galois
representation (attached to $\pi$), but this is unknown in general, except for
small $m$ (cf. \cite{Ra6, Ra5}). The difficulty here is caused by the image of Galois not
(usually) being finite.

\medskip

As mentioned in the Introduction, a result of Gelbart and Jacquet
(\cite{GeJ}) that sym$^2(\pi)$ exists for any $\pi \in \I_0(2,F)$.
It is cuspidal iff $\pi$ is not {\it dihedral}, i.e., $\pi$ is not
automorphically induced by an idele class character of a quadratic
field.

When $\pi$ is dihedral, it is easy to see that sym$^m(\pi)$ exists
for all $m$, and that it is an isobaric sum of elements of
${\mathcal A}(1,F)$ and $\I_0(2, F)$. So we may, and we will,
henceforth restrict our attention to non-dihedral forms $\pi$.

\medskip

Here is a ground-breaking result due to Kim and Shahidi which we
will need:

\medskip

\noindent{\bf Theorem 1.3.5} (Kim-Shahidi \cite{KSh1}, \cite{KSh2},
Kim \cite{K}) \, \it Let $\pi \in \I_0(2, F)$ be non-dihedral. Then
${\rm sym}^n(\pi)$ exists in ${\mathcal A}(F)$ for all $n \leq 4$.
Moreover, sym$^3(\pi)$ (resp. sym$^4(\pi)$) is cuspidal iff $\pi$ is
not {\it tetrahedral} (resp. {\it octahedral}).   \rm

\medskip

A non-dihedral $\pi$ is {\it tetrahedral} iff sym$^2(\pi)$ is
monomial, while $\pi$ is {\it octahedral} if it is not dihedral or
tetrahedral but whose symmetric cube is not cuspidal upon base
change to some quadratic extension $K$ of $F$. We will say that
$\pi$ is {\it solvable polyhedral} if it is either dihedral, or
tetrahedral, or octahedral.

\subsection{Rankin-Selberg $L$-functions}

\medskip

Let $\pi,$ $\pi'$ be isobaric automorphic representations in
$\mathcal A(n,F),$ $\mathcal A(n',F)$ respectively. Then there exists
an associated Euler product $L(s, \pi \times \pi')$ (\cite{JPSS},
\cite{JS1, JPSS, Sh1, Sh3, MW}), which converges in some right half plane,
even in $\{\Re(s)
> 1 \}$ if $\pi, \pi'$ are unitary and cuspidal. It also admits a meromorphic continuation to the whole
$s-$plane and satisfies the functional equation
$$
L(s, \pi \times \pi') \, = \, \varepsilon(s, \pi \times \pi')L(1-s,
\pi^\vee \times {\pi'}^\vee), \leqno(1.4.1)
$$
with
$$
\varepsilon(s, \pi \times \pi') \, = \, W(\pi \times \pi') N(\pi
\times \pi')^{\frac{1}{2}-s},
$$
where $N(\pi \times \pi')$ is a positive integer not divisible by
any rational prime not intersecting the ramification loci of
$F/\Q$, $\pi$ and $\pi'$, while $W(\pi \times \pi')$ is a non-zero
complex number, called the {\it root number of the pair} $(\pi,
\pi')$. As in the Galois case, $W(\pi \times \pi')W(\pi^\vee
\times {\pi'}^\vee) = 1$, so that $W(\pi \times \pi') = \pm 1$
when $\pi, \pi'$ are self-dual.

When $v$ is archimedean or a finite place unramified for $\pi,
\pi'$,
$$
L_v(s, \pi \times \pi') \, = \, L(s, \sigma(\pi_v) \otimes
\sigma(\pi'_v)). \leqno (1.4.2)
$$
In the archimedean situation, $\pi_v \to \sigma(\pi_v)$ is the
local Langlands correspondence ([La1]), with $\sigma(\pi_v)$
a representation of the Weil group $W_{F_v}$.
When $v$ is an unramified finite place, $\sigma(\pi_v)$ is
defined in the obvious way as the sum of one dimensional
representations defined by the Langlands class $A(\pi_v)$.

When $n=1,$ $L(s, \pi \times \pi') \, = \, L(s, \pi\pi'),$ and
when $n=2$ and $F = \Q,$ this function is the usual Rankin-Selberg
$L-$function, extended to arbitrary global fields by Jacquet.

\medskip

\noindent{\bf Theorem 1.4.3} \cite{JS1, JPSS}) \quad \it Let $\pi
\in \mathcal A_0(n,F)$, $\pi' \in \mathcal A_0(n',F)$, and $S$ a
finite set of places. Then $L^S(s, \pi \times \pi')$ is entire
unless $\pi$ is of the form ${\pi'}^\vee \otimes |.|^w$, in which
case it is holomorphic outside $s = -w, 1-w$, where it has simple
poles. In particular, if $\pi, \pi'$ are unitary cuspidal representations,
$L^S(s,\pi \times \pi')$ is holomorphic in $\Re(s)>1$, and moreover,
there is a pole at $s=1$ iff $\pi'\simeq\pi^\vee\simeq{\overline \pi}$.
\rm

\medskip

One also knows (cf.\cite{Sh1}, that for $\pi, \pi'$ {\it unitary} cuspidal,
$$
L^S(1+t,\pi \times \pi')\ne 0, \, \forall t\in \R.\leqno(1.4.4)
$$
Clearly, this continues to hold for isobaric sums $\pi, \pi'$ of unitary cuspidal representations.
(Note that there unitary isobaric representations which are not isobaric sums of unitary
cuspidal representations, and the assertion will not hold for these representations.)

\medskip

\subsection{The (conjectural) automorphic tensor product}

\medskip

The {\it Principle of Functoriality} asserts that given isobaric
automorphic representations $\pi, \pi'$ GL$_n(\A_F)$,
GL$_{n'}(\A_F)$ respectively, there should exist an isobaric
automorphic representation $\pi\boxtimes\pi'$, called the {\it
automorphic tensor product}, or the {\it functorial product}, of
GL$(nn',\A_F)$ such that
$$
L(s,\pi\boxtimes\pi') \, = \, L(s,\pi\times\pi').\leqno(1.5.1)
$$
We will say that an automorphic $\pi\boxtimes \pi'$ is a {\it weak
automorphic tensor product} of $\pi, \pi'$ if the identity (1.5.1)
of Euler products over $F$ holds outside a finite set $S$ of places, i.e, iff
$L(s,\pi_v\boxtimes\pi'_v)$ equals $L(s,\pi_v\times\pi'_v)$ at every $v \notin S$.

\medskip

The (conjectural) functorial product $\boxtimes$ is the automorphic
analogue of the usual tensor product of Galois representations. For
the importance of this product, see \cite{Ra1}, for example.

One can always define $\pi\boxtimes\pi'$ as an {\it admissible}
representation of GL$_{nn'}(\A_F)$, but the subtlety lies in showing
that this product is automorphic. Also, if one knows how to
construct it for cuspidal $\pi, \pi'$, then one can do it in
general.

The automorphy of $\boxtimes$ is known in the following cases, which
will be useful to us:

\noindent{${\bf (1.5.2)}$}
\begin{enumerate}
\item[]\quad $(n,n')= (2,2)$: \, (\cite{Ra2})
\item[]\quad $(n,n')=(2,3)$: \, Kim-Shahidi (\cite{KSh1})
\end{enumerate}

\medskip

The reader is also referred to section $11$ of \cite{Ra4}, which
contains some refinements, explanations, and (minor) errata for
\cite{Ra2}. Furthermore, it may be worthwhile remarking that Kim and
Shahidi effectively use their construction of the functorial product
on GL$(2)\times $GL$(3)$ to prove the automorphy of {\it symmetric
cube} transfer from GL$(2)$ to GL$(4)$, mentioned in section 1.3. A
{\it cuspidality criterion} for the image under this transfer is
proved in \cite{RaW}, with an application to the cuspidal cohomology
of congruence subgroups of SL$(6,\Z)$.

\vskip 0.2in

\section{Statement of the Main Result}

\bigskip

Here is a more precise, though a bit more cumbersome, version of Theorem
A, which was stated in the Introduction.

\medskip

\noindent{\bf Theorem A$'$} \, \it Let $\pi$ be a cuspidal automorphic
representation of GL$_2(\A_F)$ of central character $\omega$.
Assume, for the first three parts that $\pi$ is not solvable
polyhedral. Then we have the following:
\begin{enumerate}
\item[(a)] If sym$^5(\pi)$ is modular, then it is cuspidal.
\item[(b)] If sym$^5(\pi)$ and sym$^6(\pi)$ are both modular,
then sym$^6(\pi)$ is non-cuspidal iff we have
$$
{\rm sym}^5(\pi) \, \simeq \, {\rm Ad}(\pi')\boxtimes
\pi\otimes\omega^2,
$$
for a cuspidal automorphic representation $\pi'$ of GL$_2(\A_F)$; in
this case, ${\rm Ad}(\pi')$ and ${\rm Ad}(\pi)$ are not twist
equivalent.
\item[(c)] Let $m\geq 6$, and
\begin{enumerate}
\item[]assume that either
\item[(i)]sym$^j(\pi)$ is modular for every $j\leq 2m$,
\item[]or
\item[(ii)]$\pi\boxtimes \tau$ is modular for
any cusp form $\tau$ on GL$(r)/F$, with $r\leq
\left[\frac{m}{2}+1\right]$.
\end{enumerate}
Then sym$^m(\pi)$ is cuspidal iff
sym$^6(\pi)$ is cuspidal.
\item[(d)] If $F=\Q$ and $\pi$ is
defined by a non-CM, holomorphic newform $\varphi$ of weight $k\geq
2$, then ${\rm sym}^m(\pi)$ is cuspidal whenever it is modular.
\end{enumerate}
\rm

\vskip 0.2in

\section{\bf Proof of Theorem A$'$, parts (a)--(c)}

\bigskip

\subsection{Two lemmas}

In this and the following sections, $S$ will always denote a finite
set of places of $F$ containing the archimedean and finite ramified
(for $\pi$) places of $F$.

\medskip

\noindent{\bf Lemma 3.1.1} \, \it If sym$^m(\pi)$ is weakly modular, then it must
be an isobaric sum of unitary cuspidal representations. \rm

\medskip

\noindent{\it Proof}. \, Assume sym$^m(\pi)$ is weakly modular, i.e., for all places $v$
outside a finite set $S$, sym$^m(\pi_v)$ is the $v$-component of an isobaric automorphic
representation $\Pi$. Suppose $\Pi$ admits as an isobaric
summand $\Pi_0$, which is cuspidal but not unitary. In other words, there is a
non-zero real number $t$ such that $\Pi_0\otimes(\vert\cdot\vert\circ {\rm det})$
is a unitary cuspidal representation. Then {\it every} local component $\Pi_{0,v}$ is
necessarily non-unitary. As $\Pi_{0,v}$ must be a local isobaric summand of
sym$^m(\pi_v)$ for $v\notin S$,
the latter must be non-tempered.

On the other hand, since $\pi$ is a cusp form on GL$(2)/F$, we know (cf. \cite{Ra0})
that it contains infinitely many components $\pi_v$ which are tempered. (In fact, more than
$\frac{9}{10}$-th of the components are tempered.) This implies that for any finite set $S$ of places
of $F$, there exist places $v \notin S$ such that sym$^m(\pi_v)$ is tempered. This gives the
desired contradiction, yielding the Lemma.

\qed

\medskip

We will use Lemma 3.1.1 repeatedly, often without specifically referring to it.

\bigskip

\noindent{\bf Lemma 3.1.2} \, \it Suppose sym$^r(\pi)$ is modular for
all $r< m$. Pick any positive integer $i\leq m$. Then sym$^m(\pi)$
is modular iff sym$^i(\pi)\boxtimes$sym$^{m-i}(\pi)$ is modular.\rm

\medskip

{\it Proof}. \, Since $\boxtimes$ is commutative, we may assume that
$i \leq m/2$. By the Clebsch-Gordon identities, if $r_0$ denotes the
standard $2$-dimensional representation of GL$(2,\C)$, we have
$$
{\rm sym}^i(r_0)\times {\rm sym}^{m-i}(r_0) \, \simeq \,
\oplus_{j=0}^i \, {\rm sym}^{m-2j}(r_0)\otimes {\rm det}^j.
$$
It follows that
$$
L^S(s, {\rm sym}^i(\pi)\times {\rm sym}^{m-i}(\pi)) \, = \,
\prod\limits_{j=0}^i \, L^S(s, {\rm sym}^{m-2j}(\pi)\otimes
\omega^j).\leqno(3.1.3)
$$
By hypothesis, sym$^j(\pi)$ is modular for all $j<m$. If
sym$^m(\pi)$ is also modular, we may set
$$
{\rm sym}^i(\pi)\boxtimes {\rm sym}^{m-i}(\pi) := \,
\boxplus_{j=0}^i \, {\rm sym}^{m-2j}(\pi)\otimes\omega^j,
$$
which defines the desired automorphic form on GL$((i+1)(m-i+1))/F$.
Conversely,
if ${\rm sym}^i(\pi)\boxtimes {\rm sym}^{m-i}(\pi)$ is modular, then
by (3.1.3), it must have a unique isobaric summand $\Pi$, with
$$
{\rm sym}^i(\pi)\boxtimes {\rm sym}^{m-i}(\pi) := \, \Pi \boxplus
\left(\boxplus_{j=1}^i \, {\rm
sym}^{m-2j}(\pi)\otimes\omega^j\right).
$$
It follows that at any unramified place $v$ one has, for every integer $k\leq
m$ and for every irreducible admissible representation $\eta$ of
GL$_k(F_v)$, identities of the Rankin-Selberg local factors:
$$
L(s,\Pi_v\times \eta) \, = \, L(s, {\rm sym}^m(\pi)\times\eta),
$$
and
$$
\varepsilon(s,\Pi_v\times \eta) \, = \, \varepsilon(s, {\rm
sym}^m(\pi_v)\times\eta).
$$
One gets the weak modularity of sym$^m(\pi)$.
In fact, these identities hold at every place $v$, as seen by using the
local Langlands correspondence
for GL$(n)$ (\cite{HaT}, \cite{He}). From the local converse theorem, one
gets an isomorphism of $\Pi_v$
with sym$^m(\pi_v)$. Hence ${\rm sym}^m(\pi)$ is modular.

\qed

\medskip

\subsection{Proof of part (a) of Theorem A$'$}
By the work of Kim
and Shahidi (see section 1), we know that for all $j \leq 4$,
sym$^j(\pi)$ is modular, even cuspidal since $\pi$ is not solvable
polyhedral. By hypothesis, sym$^5(\pi)$ is modular. Applying Lemma
3.1.2 above with $i=4$, we get the modularity of sym$^4(\pi)\boxtimes
\pi$. Suppose sym$^5(\pi)$ is Eisensteinian. Then it must have an
isobaric summand $\tau$, say, which is cuspidal on GL$(r)/F$ for
some $r\leq 3$. By Lemma 3.1.1, $\tau$ must be unitary.
We also know (see section 1) that $\pi\boxtimes
\tau^\vee$ is automorphic on GL$(2r)/F$. Using (3.1.3) we get the
identity
$$
L^S(s,{\rm sym}^4(\pi)\times(\pi\boxtimes\tau^\vee)) \, = \,
L^S(s,{\rm sym}^5(\pi)\times\tau^\vee)L^S(s, {\rm
sym}^3(\pi)\otimes\omega\times\tau^\vee).
$$
As $\tau$ is a (unitary) cuspidal isobaric summand of sym$^5(\pi)$, the first
$L$-function on the right has a pole at $s=1$. And by the
Rankin-Selberg theory (see (1.4.4)), the second $L$-function on the right has no
zero at $s=1$. It follows that
$$
-{\rm ord}_{s=1} L^S(s,{\rm
sym}^4(\pi)\times(\pi\boxtimes\tau^\vee)) \, \geq \, 1.
$$
Since sym$^4(\pi)$ is a cusp form on GL$(5)/F$, we are forced to
have $2r\geq 5$, so $r=3$. Comparing degrees, we must then have an isobaric sum
decomposition
$$
\pi^\vee\boxtimes\tau \, \simeq \, {\rm sym}^4(\pi) \boxplus \nu,
$$
where $\nu$ is an idele class character of $F$. This implies that
$$
-{\rm ord}_{s=1} L^S(s,\pi^\vee\boxtimes\tau\otimes\nu^{-1}) \, \geq
\, 1,
$$
which is impossible unless $r=2$ and $\tau\simeq\pi\otimes\nu$. But
we have $r=3$, furnishing the desired contradiction. Hence
sym$^5(\pi)$ must be cuspidal.

\qed

\medskip

\subsection{Proof of part (b) of Theorem A$'$}
By hypothesis, sym$^j(\pi)$ is modular for all $j\leq 6$, even
cuspidal for $j\leq 5$ by part (a). By Lemma 1, ${\rm
sym}^j(\pi)\boxtimes\pi$ is also modular for each $j \leq 5$.

First suppose we have an isomorphism
$$
{\rm sym}^5(\pi)\simeq {\rm sym}^2(\pi')\boxtimes\pi\otimes\nu,
$$
for a cusp form $\pi'$ on GL$(2)/F$ and an idele class character $\nu$
of $F$. This results in the identity:
$$
L^S(s,{\rm sym}^5(\pi)\boxtimes\pi) \, = \, L^S(s,\left({\rm
sym}^2(\pi')\boxtimes\pi\right)\times\pi\otimes\nu).\leqno(3.3.1)
$$
The $L$-function on the right is the same as
$$
L^S(s,{\rm sym}^2(\pi')\times{\rm sym}^2(\pi)\otimes\nu)L^S(s,{\rm
sym}^2(\pi')\otimes\omega\nu).\leqno(3.3.2)
$$
As sym$^2(\pi')^\vee\otimes(\omega\nu)^{-1}$ is equivalent to
sym$^2(\pi')\otimes\omega\nu^{-1}$, we see that by Lemma 3.1.2,
$\Pi':={\rm sym}^2(\pi')\boxtimes{\rm
sym}^2(\pi')^\vee\otimes(\omega\nu)^{-1}$ makes sense as an
automorphic form on GL$(6)/F$. In addition, since sym$^5(\pi)\boxtimes\pi$ is
isomorphic to sym$^6(\pi)\boxplus \left({\rm
sym}^4(\pi)\otimes\omega\right)$, we obtain by using (3.3.1) and
(3.3.2):
$$
L^S(s,{\rm sym}^6(\pi)\times {\rm
sym}^2(\pi')^\vee\otimes (\omega\nu)^{-1})L^S(s,{\rm sym}^4(\pi)\times{\rm
sym}^2(\pi')^\vee\otimes\nu^{-1})\leqno(3.3.3-a)
$$
equals
$$
L^S(s,\Pi'\times{\rm sym}^2(\pi'))L^S(s,{\rm
sym}^2(\pi')\boxtimes{\rm sym}^2(\pi')^\vee).\leqno(3.3.3-b)
$$
The second $L$-function of (3.3.3-b) has a pole at $s=1$. And since
sym$^4(\pi)$ is a cusp form on GL$(5)/F$, the second $L$-function of
(3.3.3-a) has no pole at $s=1$, and the first $L$-function of
(3.3.3-b) has no zero at $s=1$. Consequently,
$$
-{\rm ord}_{s=1} L^S(s,{\rm sym}^6(\pi)\times {\rm
sym}^2(\pi')^\vee\otimes (\omega\nu)^{-1}) \, \geq \, 1.\leqno(\ast)
$$
As sym$^2(\pi')^\vee$ is automorphic on GL$(3)/F$, $(\ast)$ cannot hold
unless sym$^6(\pi)$ is not cuspidal. We are done in this direction.

\medskip

Now let us prove the {\it converse}, by supposing that sym$^6(\pi)$
is Eisensteinian. In this case it must admit an isobaric summand
$\tau$ which is cuspidal on GL$(k)/F$ with $k\leq 3$. Since we have
$$
{\rm sym}^6(\pi) \boxplus {\rm sym}^4(\pi) \, \simeq \, {\rm
sym}^5(\pi)\boxtimes\pi,
$$
$\tau$ must be an isobaric summand of sym$^5(\pi)\boxtimes\pi$. It
follows that
$$
-{\rm ord}_{s=1} L^S(s, {\rm
sym}^5(\pi)\times\left(\pi\boxtimes\tau^\vee\right)) \, \geq \, 1,
$$
where $\pi\boxtimes\tau^\vee$ is modular since $k\leq 3$. Since
sym$^5(\pi)$ is a cusp form on GL$(6)/F$, we are forced to have
$k=3$, and moreover,
$$
{\rm sym}^5(\pi) \, \simeq \, \pi^\vee\boxtimes\tau.\leqno(3.3.4)
$$
As sym$^6(\pi)$ cannot have a GL$(1)$ isobaric summand, no twist of
$\tau$ can be an isobaric summand of sym$^6(\pi)$ either, which has
degree $7$. On the other hand, since
the dual of sym$^6(\pi)$ is its twist by $\omega^{-6}$, $\tau^\vee$
is an isobaric summand of sym$^6(\pi)\otimes\omega^{-6}$. So we must
have
$$
\tau^\vee \, \simeq \, \tau\otimes\omega^{-6},\leqno(3.3.5)
$$
showing $\tau$ is essentially selfdual. In fact, if we put
$$
\eta: = \, \tau\otimes\omega^{-3},\leqno(3.3.6)
$$
it is immediate that $\eta$ is even {\it selfdual}. It follows that
$$
L^S(s,\eta,{\rm sym}^2)L^S(s,\eta,
\Lambda^2)=L^S(s,\eta\times\eta^\vee),
$$
showing that the left hand side has a pole at $s=1$. Since $\eta$ is
a cusp form on GL$(3)/F$, the second $L$-function cannot have a pole
at $s=1$ (see \cite{JS2}). Hence
$$
-{\rm ord}_{s=1} L^S(s, \eta,{\rm sym}^2) \, \geq \, 1.\leqno(3.3.7)
$$
By the backwards lifting results of Ginzburg, Rallis and Soudry
(cf. \cite{GRS}, \cite{Sou}), we then have a functorially associated cuspidal,
necessarily generic, automorphic representation $\pi'_0$ of
SL$(2,\A_F)$ ($={\rm Sp}(2,\A_F)$) of trivial central character. We
may extend it (cf. \cite{LL}) to an irreducible cusp
form $\pi'$ of GL$(2)/F$ (with central character $\omega'$), which is unique only
up to twisting by a character, such that
$$
L^S(s, {\rm Ad}(\pi')) \, = \, L^S(s, \eta).\leqno(3.3.7)
$$
By the strong multiplicity one theorem, $\eta$ is isomorphic to
Ad$(\pi')$, which is sym$^2(\pi')\otimes\omega'^{-1}$.

Combining with (3.3.4) and (3.3.6), we get
$$
{\rm sym}^5(\pi) \, \simeq \, {\rm Ad}(\pi')\boxtimes
\pi\otimes\omega^{2},
$$
as asserted in part (b) of Theorem A$'$.

Finally suppose sym$(\pi)$ and Ad$(\pi')$ are twist equivalent. Then
sym$^5(\pi)$ would need to be twist equivalent to
sym$^2(\pi)\boxtimes \pi$, which is Eisensteinian of the form
sym$^3(\pi)\boxplus \pi\otimes\omega$. This contradicts the
cuspidality of sym$^5(\pi)$, and we are done.

\qed

\bigskip

\subsection{Proof of part (c) under Assumption (i)}

\medskip

There is nothing to prove if $m=6$, so let $m\geq 7$, and assume by
induction that the conclusion holds for all $n\leq m-1$. In
particular, sym$^n(\pi)$ is cuspidal for every $n<m$. Moreover, by
hypothesis, sym$^j(\pi)$ is modular for all $j \leq 2m$, and this
implies, by Lemma 3.1.2, that sym$^m(\pi)\boxtimes{\rm sym}^m(\pi)$ is
modular.

Suppose sym$^m(\pi)$ is not cuspidal. Then by \cite{JS1},
$$
-{\rm ord}_{s=1} L^S(s, {\rm sym}^m(\pi)\times {\rm
sym}^m(\pi)^\vee) \, \geq \, 2.\leqno(3.4.1)
$$
We have by Clebsch-Gordon,
$$
{\rm sym}^m(\pi)\boxtimes {\rm sym}^m(\pi)^\vee \, \simeq \,
\boxplus_{j=0}^m \, {\rm sym}^{2j}(\pi)\otimes\omega^{-j},
$$
and of course we have a similar formula for ${\rm
sym}^{m-1}(\pi)\boxtimes {\rm sym}^{m-1}(\pi)^\vee$, where the sum
goes from $j=0$ to $j=m-1$. Consequently,
$$
{\rm sym}^m(\pi)\boxtimes {\rm sym}^m(\pi)^\vee \, \simeq \,
\left({\rm sym}^{m-1}(\pi)\boxtimes {\rm sym}^{m-1}(\pi)^\vee\right)
\boxplus \left({\rm
sym}^{2m}(\pi)\otimes\omega^{-m}\right).\leqno(3.4.2)
$$
Since sym$^{m-1}(\pi)$ is cuspidal, $L^S(s, {\rm
sym}^{m-1}(\pi)\times {\rm sym}^{m-1}(\pi)^\vee)$ has a simple pole
at $s=1$ (cf. \cite{JS1}). Combining this with (3.4.1) and (3.4.2),
we obtain
$$
-{\rm ord}_{s=1} L^S(s, {\rm sym}^{2m}(\pi)\otimes\omega^{-m}) \,
\geq \, 1.\leqno(3.4.3)
$$
Since ${\rm sym}^{2m}(\pi)$ is automorphic, it must admit $\omega^m$
as an isobaric summand.

On the other hand, we have (by Clebsch-Gordon)
$$
{\rm sym}^{m+1}(\pi)\boxtimes{\rm sym}^{m-1}(\pi) \, \simeq \,
\boxplus_{j=0}^{m-1} \, {\rm
sym}^{2(m-j)}(\pi)\otimes\omega^{j}.\leqno(3.4.4)
$$
It follows that $\omega^m$ must be an isobaric summand of ${\rm
sym}^{m+1}(\pi)\boxtimes{\rm sym}^{m-1}(\pi)$, implying
$$
-{\rm ord}_{s=1} L^S(s, {\rm sym}^{m+1}(\pi)\times\left({\rm
sym}^{m-1}(\pi)\otimes\omega^{-m}\right)) \, \geq \, 1.\leqno(3.4.5)
$$
Since ${\rm sym}^{m-1}(\pi)$ is cuspidal, this can only happen (cf.
\cite{JS1}) if ${\rm sym}^{m-1}(\pi)^\vee\otimes\omega^m$ is an
isobaric summand of ${\rm sym}^{m+1}(\pi)$. Therefore
$$
{\rm sym}^{m+1}(\pi) \, \simeq \, \left({\rm
sym}^{m-1}(\pi)^\vee\otimes\omega^m\right) \boxplus \tau,
$$
where $\tau$ is an (isobaric) automorphic form on GL$(2)/F$.

Hence $\tau$ is an isobaric summand of ${\rm
sym}^m(\pi)\boxtimes{\pi}$, which is isomorphic to
sym$^{m+1}(\pi)\boxplus\left({\rm
sym}^{m-1}(\pi)\otimes\omega\right)$. Recall that $\pi^\vee\boxtimes
\tau$ is modular.  Then there is an isobaric summand $\beta$ of
$\pi^\vee\boxtimes\tau$, which is cuspidal on GL$(r)/F$ with $r\leq
4$, such that
$$
-{\rm ord}_{s=1} L^S(s, {\rm sym}^m(\pi)\times\beta^\vee) \, \geq \,
1.
$$
In other words, $\beta$ is an isobaric summand of ${\rm
sym}^m(\pi)$, and hence of sym$^{m-1}(\pi)\boxtimes\pi$.
Consequently,
$$
-{\rm ord}_{s=1} L^S(s, \left({\rm
sym}^{m-1}(\pi)\boxtimes\pi\right)\times\beta^\vee) \, \geq \,
1.\leqno(3.4.6)
$$

First suppose $r\leq 3$. Then we know that $\pi\boxtimes\beta^\vee$
is modular on GL$(2r)$ (by \cite{Ra2} for r=2, and \cite{KSh1} for
$r=3$). As sym$^{m-1}(\pi)$ is by induction cuspidal, (3.4.6) forces
the bound
$$
m \, \leq \, 2r \, \leq 6.\leqno(3.4.7)
$$
So we are done in this case.

Next suppose that $r=4$, which means $\beta=\pi^\vee\boxtimes\tau$
is cuspidal. Since $\pi\boxtimes\pi^\vee \simeq {\rm
sym}^2(\pi)\boxplus \omega$, it follows that
$\pi\boxtimes\beta^\vee$ is modular, with
$$
\pi\boxtimes\beta^\vee \, \simeq \, \left({\rm
sym}^2(\pi)\boxtimes\tau^\vee\right) \boxplus
\left(\omega\otimes\tau^\vee\right),
$$
where the first summand is on GL$(6)/F$ and the second on GL$(4)$.
As a result, we have from (3.4.6),
$$
-{\rm ord}_{s=1} L^S(s, {\rm sym}^{m-1}(\pi)\times\delta) \, \geq \,
1,\leqno(3.4.8)
$$
for an isobaric summand $\delta$ of $\pi\boxtimes\beta^\vee$, which
is a cusp form on GL$(n)$, for some $n\leq 6$. So, once again, the
inequality (3.4.7) holds, and we are done.

\qed

\bigskip

\subsection{Proof of part (c) under Assumption (ii)}

\medskip

The proof of part (c) in this case is a bit different because we are
not assuming good properties of sym$^j(\pi)$ for
$j$ all the way up to $2m$.

We may take $m>6$ and assume by induction that sym$^j(\pi)$ is
cuspidal for all $j \leq m-1$. Suppose sym$^m(\pi)$ is
Eisensteinian. Then it must have an isobaric summand $\eta$, which
is cuspidal on GL$(r)/F$ with $r\leq \left[\frac{m+1}{2}\right]$. Then $\eta$ must be an
isobaric summand of sym$^{m-1}(\pi)\boxtimes\pi$, because of the
decomposition
$$
{\rm sym}^{m-1}(\pi)\boxtimes\pi \, \simeq \, {\rm sym}^m(\pi)
\boxplus \left({\rm sym}^{m-2}(\pi)\otimes\omega\right).
$$
By our hypothesis, $\pi\boxtimes \eta^\vee$ is modular on
GL$(2r)/F$. So we get
$$
-{\rm ord}_{s=1} L^S(s, {\rm sym}^{m-1}(\pi)\times\left(\pi\boxtimes
\eta^\vee\right)) \, \geq 1.\leqno(3.5.1)
$$
As sym$^{m-1}(\pi)$ is cuspidal, we are forced to have
$$
m \, \leq \, 2r \, \leq \, m+1.\leqno(3.5.2)
$$
So the only possible (isobaric) decomposition of sym$^m(\pi)$ we can
have is
$$
{\rm sym}^m(\pi) \, \simeq \, \eta \boxplus \eta', \leqno(3.5.3)
$$
with
$$
\eta \in \I_0([(m+1)/2], F) \quad {\rm and} \quad \eta' \in
\I_0(m+1-[(m+1)/2], F).
$$
In addition, by our hypothesis, $\eta \boxtimes \pi^\vee$ and $\eta'
\boxtimes \pi^\vee$ are modular. We deduce that
$$
[{\rm sym}^{m-1}(\pi), \eta \boxtimes \pi^\vee] > 0, \quad {\rm and}
\quad [{\rm sym}^{m-1}(\pi), \eta' \boxtimes \pi^\vee] > 0.
\leqno(3.5.4)
$$

First consider the case when {\it $m$ is odd}. (This is similar to
the argument above for $m = 5$.) Then $r = [(m+1)/2] = m+1
-[(m+1)/2]$, and $\eta, \eta'$ are both in $\I_0((m+1)/2, F)$.
Since sym$^{m-1}(\pi) \in \I_0(m, F)$, we must have
$$
\eta \boxtimes \pi^\vee \, \simeq \, {\rm sym}^{m-1}(\pi) \boxplus
\mu
$$
and
$$
\eta' \boxtimes \pi^\vee \, \simeq \, {\rm sym}^{m-1}(\pi) \boxplus
\mu',
$$
with $\mu, \mu'$ in ${\mathcal A}(1,F)$. Then it follows that the
Rankin-Selberg $L$-functions $L^S(s, \eta \times (\pi^\vee \otimes
\mu^{-1}))$ and $L^S(s, \eta' \times (\pi^\vee \otimes
{\mu'}^{-1}))$ both have poles at $s=1$. This forces the following:
$$
m=3, \, \, \eta \simeq \pi \otimes \mu, \quad {\rm and} \quad \eta'
\simeq \pi \otimes \mu'.
$$
So this cannot happen for $m \ne 3$.

Next consider the case when {\it $m$ is even}. Then $\eta \in
\I_0(m/2, F)$ and $\eta' \in \I_0(m/2 + 1, F)$. We get
$$
\eta \boxtimes \pi^\vee \, \simeq \, {\rm sym}^{m-1}(\pi)
$$
and
$$
\eta' \boxtimes \pi^\vee \, \simeq \, {\rm sym}^{m-1}(\pi) \boxplus
\tau,
$$
with $\tau$ in $\I_0(2,F)$. Then $\eta'$ must occur in $\pi
\boxtimes \tau$, which is in ${\mathcal A}(4, F)$. So we must have
$$
m/2 + 1 \, \leq \, 4.
$$
In other words, $m$ must be less than or equal to $6$, which is not
the case.

Thus we get a contradiction in either case. The only possibility is
for sym$^m(\pi)$ to be cuspidal. Done proving part (c), and hence
all of Theorem B.

\qed

\bigskip

\section{Proof of Theorem A$'$, part (d)}

\medskip

Finally, we want to restrict to $F = \Q$ and analyze the case of
holomorphic newforms $f$ of weight $\geq 2$. One knows that the
level $N$ of $f$ is the same as the conductor of the associated
cuspidal automorphic representation $\pi$ of GL$(2,\A_\Q)$.
Moreover, as $f$ is not of CM type, $\pi$ is not dihedral.

\medskip

Fix a prime $\ell$ not dividing $N$ and consider the
{\it cyclotomic character}
$$
\chi_\ell: {\text Gal}(\overline \Q/\Q) \, \rightarrow \,
\Z_\ell^\ast,\leqno(4.1)
$$
defined by the Galois action on the projective system
$\{\mu_{\ell^r} \vert r \geq 1\}$, where $\mu_{\ell^r}$ denotes the
group of $\ell^r$-th roots of unity in $\overline \Q$. Then by a
theorem of Deligne, one has at our disposal an irreducible,
continuous representation
$$
\rho_\ell(\pi): {\text Gal}(\overline \Q/\Q) \, \rightarrow \,
{\text GL}(2, \overline \Q_\ell), \leqno(4.2)
$$
unramified outside $N\ell$,
such that for every prime $p$ not dividing $N\ell$,
$$
{\text Tr}(\rho_\ell(\pi)({\text Fr}_p)) \, = \, a_p,\leqno(4.3)
$$
where Fr$_p$ denotes the Frobenius at $p$ and $a_p$ the $p$-th Hecke
eigenvalue of $f$. Moreover,
$$
{\text det}(\rho_\ell(\pi) \, = \, \omega\chi_\ell^{k-1}.\leqno(4.4)
$$
When $f$ is of CM-type, there exists an imaginary quadratic field
$K$, and an algebraic Hecke character $\Psi$ of $K$ such that
$$
\rho_\ell(\pi) \, \simeq \, {\text Ind}_{{\text Gal}(\overline
\Q/K)}^ {{\text Gal}(\overline \Q/\Q)} (\Psi_\ell),\leqno(4.5)
$$
where $\Psi_\ell$ is the $\ell$-adic character associated to $\Psi$
(\cite{Se}). Let $\theta$ denote the non-trivial automorphism of
Gal$(K/\Q)$. Then it is an immediate exercise to check that for any
$m \geq 1$, sym$^m(\rho_\ell)$ is of the form $\oplus_j
\beta_{j,\ell}$, where each $\beta_{j,\ell}$ is either
one-dimensional defined by an idele class character of $\Q$ or a
two-dimensional induced by $\Psi_\ell^{a}(\Psi_\ell^\theta)^{m-a}$
for some $a \geq 0$, with $\Psi_\ell^\theta$ denoting the conjugate
of $\Psi_\ell$ under $\theta$. It is clear this is modular, but not
cuspidal for any $m\geq 2$.

Let us assume henceforth that $f$ is not of $CM$-type. Denote by
$G_\ell$ the {\it Zariski closure of the image} of Gal$(\overline
\Q/\Q)$ under $\rho_\ell(\pi)$; it is an $\ell$-adic Lie group.
Since $f$ is of weight $\geq 2$ and {\it not} of CM-type, a theorem
of K.~Ribet (\cite{Ri}) asserts that for {\it large enough} $\ell$,
$$
G_\ell \, = \, {\text GL}(2, \overline \Q_\ell). \leqno(4.6)
$$
We will from now on consider only those $\ell$ large enough for this
to hold. Since the symmetric power representations of the algebraic
group GL$(2)$ are irreducible, we get the following

\medskip

\noindent{\bf Lemma 4.7} \, \it For any non-CM newform $f$ of weight
$k \geq 2$ and for any $m \geq 1$ and large enough $\ell$, the
representation $\s^m(\rho_\ell)$ is irreducible, and it remains so
under restriction to Gal$(\overline \Q/E)$ for any finite extension
$E$ of $\Q$. \rm

\medskip

Since $f$ is not of CM-type, $\s^2(\pi)$ is cuspidal. In view of
parts (a)--(c) (of Theorem A$'$), we need only prove the following
to deduce part (d):

\medskip

\noindent{\bf Proposition 4.8} \, \it For any non-CM newform $f$ of
weight $k \geq 2$ and level $N$, with associated cuspidal
automorphic representation $\pi$ of GL$(2,\A_\Q)$, assume that
sym$^m(\pi)$ is modular for all $m\geq 2$. Then the following hold:
\begin{enumerate}
\item[(i)]{For any quadratic field $K$, the base change $\s^3(\pi)_K$
to GL$(4)/K$ is cuspidal}
\item[(ii)]{$\s^6(\pi)$ is cuspidal}
\end{enumerate}
\rm

\medskip

This Proposition suffices, because (i) implies that $\pi$ is not
solvable polyhedral, and (ii) implies what we want by part (c) of
Theorem A$'$.

\medskip

Let $f$ be as in the Proposition. Suppose $m \geq 1$ is such that
$\s^j(\pi)$ is cuspidal for all $j < m$, but Eisensteinian for
$j=m$. Then we have, as in the proof of the earlier parts of Theorem
A$'$, a decomposition
$$
{\rm sym}^m(\pi) \, \simeq \, \eta \boxplus \eta', \leqno(4.9)
$$
with
$$
\eta \in {\mathcal A}_0([(m+1)/2], \Q) \quad {\rm and} \quad \eta'
\in {\mathcal A}_0(m+1-[(m+1)/2], \Q),
$$
with $\eta, \eta'$ are essentially self-dual. Moreover, we have

\medskip

\noindent{\bf Lemma 4.10} \, \it The infinity types of $\eta, \eta'$
are both algebraic and regular.\rm

\medskip

Some explanation of the terminology is called for at this point.
Recall that $W_\R$ is the unique non-split extension of Gal$(\C/\R)$
by $\C^\ast$, which is concretely described as $\C^\ast \cup
j\C^\ast$, with $jzj^{-1} = \overline z$, for all $z \in \C^\ast$.
Let $\Pi$ be an irreducible automorphic representation of
GL$(n,\A_F)$. Since the restriction of $\sigma_\infty(\Pi)$ is
semisimple and since $\C^\ast$ is abelian, we get a decomposition
$$
\sigma_\infty(\Pi)\vert_{\C^\ast} \, \simeq \, \oplus_{i \in J} \chi_i,
$$
where each $\chi_i$ is in Hom$_{\text cont}(\C^\ast, \C^\ast)$.
$\Pi_\infty$ is said to be {\bf regular} iff this decomposition is
multiplicity-free, i.e., iff $\chi_i \ne \chi_r$ for $i \ne r$. It
is {\bf algebraic} (\cite{C}) iff each
$\chi_i\vert\cdot\vert^{(m-1)/2}$ is of the form $z \to
z^{-a_i}{\overline z}^{-b_i}$, for some integers $a_i, b_i$. An
algebraic $\Pi$ is said to be {\bf pure} if there is an integer $w$,
called the {\bf weight} of $\Pi$, such that $w = a_i + b_i$ for each
$i \in J$.

It is well known that, since $\pi$ is defined by a holomorphic
newforms $f$ of weight $k \geq 2$,
$$
\sigma_\infty(\pi)\otimes\vert\cdot\vert^{-1/2} \, \simeq \, {\text
Ind}(W_\R, \C^\ast; z_{1-k}), \leqno(4.11)
$$
where $z_n$ denotes, for each integer $n$, the continuous
homomorphism $\C^\ast \to \C^\ast$ given by $z \to z^n$. Note that
$\pi_\infty$ is regular (as $k > 1$) and algebraic of weight $k-1$.
From here on to the end of this chapter, we will simply write $I(-)$
for ${\text Ind}(W_\R, \C^\ast; -)$. Set
$$
\nu_{1-k} \, = \, z_{1-k}\vert_{\R^\ast}.
$$
Then we have
$$
\omega_\infty \, = \, ({\text sgn})\nu_{1-k}, \leqno(4.12)
$$
where ${\text sgn}$ denotes the sign character of $\R^\ast$.
Indeed, $\omega_\infty = {\text sgn}^{1-k}\nu_{1-k}$. But as $f$ has trivial
character, $k$ is forced to be
even, so ${\text sgn}^{1-k} = {\text sgn}$. (Here we have identified,
as we may, $\omega_\infty$
with $\sigma_\infty(\omega)$.)

\medskip

\noindent{\bf SubLemma 4.13} \, \it For each $j \leq [m/2]$,
$$
\sigma_\infty({\text sym}^{2j+1}(\pi)) \, \simeq \,
I(z_{1-k}^{2j+1}) \oplus (I(z_{1-k}^{2j-1}) \otimes |.|^{1-k}) \oplus \ldots
\oplus (I(z_{1-k}) \otimes |.|^{(1-k)j}),
\leqno(i)
$$
and
$$
\sigma_\infty({\text sym}^{2j}(\pi)) \, \simeq \,
I(z_{1-k}^{2j}) \oplus (I(z_{1-k}^{2j-2}) \otimes |.|^{1-k}) \oplus \ldots
\oplus (I(z_{1-k}^2) \otimes |.|^{(1-k)(j-1)}) \oplus \nu_{1-k}^j.
\leqno(ii)
$$
\rm

\medskip

{\it Proof of SubLemma}. \, Everything is fine for $j = 0$. So we
may let $j > 0$ and assume by induction that the identities hold for
all $r < j$. Applying (i) for $j-1$ together with $(4.3)_{2j}$,
(4.11) and (3.19), we see that
$$
\sigma_\infty({\text sym}^{2j}(\pi))  \oplus
(\sigma_\infty({\text sym}^{2j-2}(\pi)) \otimes |.|^{1-k})
$$
is isomorphic to
$$
(I(z_{1-k}^{2j-1}) \oplus (I(z_{1-k}^{2j-3}) \otimes |.|^{1-k}) \oplus \ldots
\oplus (I(z_{1-k}) \otimes |.|^{(1-k)(j-1)}) \otimes I(z_{1-k}).
$$
By Mackey theory, we have for all $a \geq b$,
$$
I(z_{1-k}^a) \otimes I(z_{1-k}^b) \, \simeq \,
I(z_{1-k}^{a+b}) \oplus I(z_{1-k}^a{\overline z}_{1-k}^b) \,
\simeq \, I(z_{1-k}^{a+b}) \oplus (I(z_{1-k}^{a-b}) \otimes |.|^{1-k}).
$$
Since $I(-) \otimes {\text sgn} \simeq I(-)$, $I(-) \otimes
|.|^{1-k}$ is isomorphic to $I(-) \otimes \nu_{1-k}$. Combining
these and using the inductive assumption for $\sigma({\text
sym}^{2j-2}(\pi))$, we get (ii) for $j$. The proof of (ii) is
similar and left to the reader.

\qed

\medskip

Now Lemma 4.10 follows easily from the SubLemma and the definition
of regular algebraicity.

\medskip

\noindent{\it Proof of Proposition} (contd.) \,  We need only
examine sym$^m(\pi)$ for $m=3$ and $m=6$.

First suppose $m = 3$. Let $K$ be any quadratic field. Then $\eta_K$
and $\eta'_K$ are both essentially self-dual forms on GL$(2)/K$ with
algebraic, regular infinity types. Consequently, one knows that for
$\beta \in \{\eta, \eta'\}$, there exists a semisimple
representation
$$
\rho_\ell(\beta): \, {\text Gal}(\overline \Q/K) \, \rightarrow \,
GL(2, \overline \Q_\ell)
$$
such that for primes $P$ in a set of Dirichlet density $1$, we have
$$
L(s, \beta_P) \, = \, {\text det}(1 - {\text Fr}_P (NP)^{-s}\vert
\rho_\ell(\beta))^{-1}. \leqno(4.14)
$$
If $\beta$ is Eisensteinian, which in fact cannot happen, this is
easy to establish. Ditto if it is dihedral. So we may take $\beta$
to be cuspidal and non-dihedral. If $K$ is totally real, the
existence of  $\rho_\ell(\beta)$ is a well known result, due
independently to R.~Taylor (\cite{Ta1}) and to Blasius-Rogawski
(\cite{BRog}); in fact a stronger assertion holds in that case. In
this case, $\beta$ corresponds to a Hilbert modular form, either one
of weight $3k-2$ or to a twist of one of weight $3k-4$. If $K$ is
imaginary, the existence of $\rho_\ell(\beta)$ is a theorem of
R.~Taylor (\cite{Ta2}), partly based on his joint work with
M.~Harris and D.~Soudry. (Note that here, the central character of
the unitary version of $\beta$ is trivial.)

By part (a) of the Lemma, we then get the following at all primes $P$
in a set of density $1$:
$$
L(s, {\text sym}^3(\pi_K)_P) \, = \, {\text det} (1 - {\text Fr}_P
(NP)^{-s}\vert \rho_\ell(\eta) \oplus \rho_\ell(\eta'))^{-1}.
\leqno(4.15)
$$
But by construction,
$$
L(s, {\text sym}^3(\pi_K)_P) \, = \, {\text det} (1 - {\text Fr}_P
(NP)^{-s}\vert {\text sym}^3(\rho_\ell(\pi)_K))^{-1}. \leqno(4.16)
$$
Thus we have,
by the Tchebotarev density theorem,
$$
{\text sym}^3(\rho_\ell(\pi)_K) \, \simeq \, \rho_\ell(\eta)
\oplus \rho_\ell(\eta').
$$
We get a contradiction as we know (cf. Lemma 4.7) that ${\text
sym}^3(\rho_\ell(\pi)_K)$ is an irreducible representation.

Thus sym$^3(\pi_K)$ is cuspidal. This proves part (i) of the
Proposition, and implies that $\pi$ is not solvable polyhedral.

\medskip

Next we turn to the question of cuspidality of $\s^6(\pi)$. Again,
thanks to the hypothesis of modularity $\s^6(\pi)$, $\s^j(\pi)$ is
cuspidal for all $j \leq 5$.

Suppose $\s^6(\pi)$ is not cuspidal. Let $\eta, \eta'$ be as in the
the decomposition $\s^m(\pi)$ given by (4.9). Since $m=6$, $\eta \in
{\mathcal A}_0(3, \Q)$ and $\eta' \in {\mathcal A}_0(4, \Q)$.
Specializing Lemma 3.1.2 to $(i,m)=(5,6)$, we get
$$
\s^5(\pi) \boxtimes \pi \, \simeq \, \eta \boxplus \eta' \boxplus
(\s^4(\pi) \otimes |.|^{1-k}). \leqno(4.17)
$$

\medskip

\noindent{\bf Lemma 4.18} \, Let $\beta \in \{\eta, \eta'\}$. Take
$m = 3$ if $\beta = \eta$ and $m = 4$ if $\beta = \eta'$.  Then for
any prime $\ell$ away from the ramification locus of $\beta$, there
exists a semisimple $\ell$-adic representation
$$
\rho_\ell(\beta): {\text Gal}(\overline \Q/\Q) \, \rightarrow \,
{\text GL}(m, \overline \Q_\ell)
$$
such that for almost all primes $p$, we have
$$
L(s, \beta_p) \, = \, {\text det}(1 - {\text Fr}_p p^{-s}\vert
\rho_\ell(\beta))^{-1}. \leqno(4.19)
$$

\medskip

{\it Proof of Lemma}. \, First Note that since the dual of
sym$^6(\pi)$ is sym$^6(\pi)\otimes\omega^{-6}$, the twisted
representation sym$^6(\pi)\otimes\omega^{-3}$ is selfdual. So, we
may, after replacing sym$^6(\pi)$, $\eta$ and $\eta'$ by their
respective twists by $\omega^{3}$, assume that they are all
selfdual. (Since $\eta, \eta'$ are irreducible representations of
unequal dimensions, they cannot be contragredients of each other,
and so are forced to be selfdual themselves.)  As we have seen, they
are also regular and algebraic. Now the discussion in \cite{Ra6}
explains how to deduce the existence of the desired Galois
representations attached to $\eta, \eta'$ (see also \cite{RaSh, Ra3,
Lau, W}).

\qed

\medskip

{\it Proof of Proposition 4.8 (contd.)}. \, Applying Lemma 4.18 we
get for almost all primes $p$,
$$
L(s, {\text sym}^6(\pi)_p) \, = \, {\text det} (1 - {\text Fr}_p
p^{-s}\vert \rho_\ell(\eta) \oplus \rho_\ell(\eta'))^{-1}.
$$
By the Tchebotarev density theorem,
$$
{\text sym}^6(\rho_\ell(\pi)) \, \simeq \, \rho_\ell(\eta) \oplus \rho_\ell(\eta').
$$
Again we get a contradiction since by Lemma 4.7, ${\text
sym}^6(\rho_\ell(\pi))$ is an irreducible representation.

Thus sym$^6(\pi)$ is cuspidal.

\qed

We have now completely proved Theorem A$'$, which implies Theorem A
of the Introduction.

\vskip 0.2in

\bibliographystyle{math}    
\bibliography{cusp-sym}

\def\cprime{$'$}
\begin{thebibliography}{JPSS2}

\bibitem[BR]{BRog}
D.~Blasius and J.~D. Rogawski.
\newblock {Motives for {H}ilbert modular forms}.
\newblock {\em Invent. Math.} {\bf 114} (1993), 55--87.

\bibitem[Clo]{C}
L.~Clozel.
\newblock {Motifs et formes automorphes: applications du principe de
  fonctorialit\'e}.
\newblock In {\em Automorphic forms, Shimura varieties, and $L$-functions,
  Vol.\ I (Ann Arbor, MI, 1988)}, volume~10 of {\em Perspect. Math.}, pages
  77--159. Academic Press, Boston, MA, 1990.

\bibitem[GJ]{GeJ}
S.~Gelbart and H.~Jacquet.
\newblock {A relation between automorphic representations of {${\rm GL}(2)$}
  and {${\rm GL}(3)$}}.
\newblock {\em Ann. Sci. \'Ecole Norm. Sup. (4)} {\bf 11} (1978), 471--542.

\bibitem[GRS]{GRS}
D.~Ginzburg, S.~Rallis, and D.~Soudry.
\newblock {On explicit lifts of cusp forms from {${\rm GL}\sb m$} to classical
  groups}.
\newblock {\em Ann. of Math. (2)} {\bf 150} (1999), 807--866.

\bibitem[GJ]{GoJ}
R.~Godement and H.~Jacquet.
\newblock {\em Zeta functions of simple algebras}.
\newblock Springer-Verlag, Berlin, 1972.
\newblock Lecture Notes in Mathematics, Vol. 260.

\bibitem[HT]{HaT}
M.~Harris and R.~Taylor.
\newblock {\em The geometry and cohomology of some simple {S}himura varieties},
  volume 151 of {\em Annals of Mathematics Studies}.
\newblock Princeton University Press, Princeton, NJ, 2001.
\newblock With an appendix by Vladimir G. Berkovich.

\bibitem[Hen]{He}
G.~Henniart.
\newblock {Une preuve simple des conjectures de {L}anglands pour {${\rm
  GL}(n)$} sur un corps {$p$}-adique}.
\newblock {\em Invent. Math.} {\bf 139} (2000), 439--455.

\bibitem[JPSS1]{JPSS2}
H.~Jacquet, I.~I. Piatetski-Shapiro, and J.~Shalika.
\newblock {Conducteur des repr\'esentations du groupe lin\'eaire}.
\newblock {\em Math. Ann.} {\bf 256} (1981), 199--214.

\bibitem[JPSS2]{JPSS}
H.~Jacquet, I.~I. Piatetskii-Shapiro, and J.~A. Shalika.
\newblock {Rankin-{S}elberg convolutions}.
\newblock {\em Amer. J. Math.} {\bf 105} (1983), 367--464.

\bibitem[JS]{JS1}
H.~Jacquet and J.~A. Shalika.
\newblock {On {E}uler products and the classification of automorphic forms.
  {II}}.
\newblock {\em Amer. J. Math.} {\bf 103} (1981), 777--815.

\bibitem[Jac]{J}
H.~Jacquet.
\newblock {Principal {$L$}-functions of the linear group}.
\newblock In {\em Automorphic forms, representations and $L$-functions (Proc.
  Sympos. Pure Math., Oregon State Univ., Corvallis, Ore., 1977), Part 2},
  Proc. Sympos. Pure Math., XXXIII, pages 63--86. Amer. Math. Soc., Providence,
  R.I., 1979.

\bibitem[JS]{JS2}
H.~Jacquet and J.~Shalika.
\newblock {Exterior square {$L$}-functions}.
\newblock In {\em Automorphic forms, Shimura varieties, and $L$-functions,
  Vol.\ II (Ann Arbor, MI, 1988)}, volume~11 of {\em Perspect. Math.}, pages
  143--226. Academic Press, Boston, MA, 1990.

\bibitem[Kim]{K}
H.~H. Kim.
\newblock {Functoriality for the exterior square of {${\rm GL}\sb 4$} and the
  symmetric fourth of {${\rm GL}\sb 2$}}.
\newblock {\em J. Amer. Math. Soc.} {\bf 16} (2003), 139--183 (electronic).
\newblock With appendix 1 by Dinakar Ramakrishnan and appendix 2 by Kim and
  Peter Sarnak.

\bibitem[KS1]{KSh2}
H.~H. Kim and F.~Shahidi.
\newblock {Cuspidality of symmetric powers with applications}.
\newblock {\em Duke Math. J.} {\bf 112} (2002), 177--197.

\bibitem[KS2]{KSh1}
H.~H. Kim and F.~Shahidi.
\newblock {Functorial products for {${\rm GL}\sb 2\times{\rm GL}\sb 3$} and the
  symmetric cube for {${\rm GL}\sb 2$}}.
\newblock {\em Ann. of Math. (2)} {\bf 155} (2002), 837--893.
\newblock With an appendix by Colin J. Bushnell and Guy Henniart.

\bibitem[Lan1]{La2}
R.~P. Langlands.
\newblock {Automorphic representations, {S}himura varieties, and motives. {E}in
  {M}\"archen}.
\newblock In {\em Automorphic forms, representations and $L$-functions (Proc.
  Sympos. Pure Math., Oregon State Univ., Corvallis, Ore., 1977), Part 2},
  Proc. Sympos. Pure Math., XXXIII, pages 205--246. Amer. Math. Soc.,
  Providence, R.I., 1979.

\bibitem[Lan2]{La3}
R.~P. Langlands.
\newblock {On the notion of an automorphic representation}.
\newblock In {\em Automorphic forms, representations and $L$-functions (Proc.
  Sympos. Pure Math., Oregon State Univ., Corvallis, Ore., 1977), Part 2},
  Proc. Sympos. Pure Math., XXXIII, pages 189--217. Amer. Math. Soc.,
  Providence, R.I., 1979.

\bibitem[Lan3]{La1}
R.~P. Langlands.
\newblock {On the classification of irreducible representations of real
  algebraic groups}.
\newblock In {\em Representation theory and harmonic analysis on semisimple Lie
  groups}, volume~31 of {\em Math. Surveys Monogr.}, pages 101--170. Amer.
  Math. Soc., Providence, RI, 1989.

\bibitem[Lan4]{La}
R.~P. Langlands.
\newblock {Beyond endoscopy}.
\newblock In {\em Contributions to automorphic forms, geometry, and number
  theory}, pages 611--697. Johns Hopkins Univ. Press, Baltimore, MD, 2004.

\bibitem[Lau]{Lau}
G.~Laumon.
\newblock {Fonctions z\^etas des vari\'et\'es de {S}iegel de dimension trois}.
\newblock {\em Ast\'erisque} (2005), 1--66.
\newblock Formes automorphes. II. Le cas du groupe $\rm GSp(4)$.

\bibitem[MW]{MW}
C.~M{\oe}glin and J.-L. Waldspurger.
\newblock {Le spectre r\'esiduel de {${\rm GL}(n)$}}.
\newblock {\em Ann. Sci. \'Ecole Norm. Sup. (4)} {\bf 22} (1989), 605--674.

\bibitem[Ram1]{Ra1}
D.~Ramakrishnan.
\newblock {Pure motives and automorphic forms}.
\newblock In {\em Motives (Seattle, WA, 1991)}, volume~55 of {\em Proc. Sympos.
  Pure Math.}, pages 411--446. Amer. Math. Soc., Providence, RI, 1994.

\bibitem[Ram2]{Ra0}
D.~Ramakrishnan.
\newblock {On the coefficients of cusp forms}.
\newblock {\em Math. Res. Lett.} {\bf 4} (1997), 295--307.

\bibitem[Ram3]{Ra2}
D.~Ramakrishnan.
\newblock {Modularity of the {R}ankin-{S}elberg {$L$}-series, and multiplicity
  one for {${\rm SL}(2)$}}.
\newblock {\em Ann. of Math. (2)} {\bf 152} (2000), 45--111.

\bibitem[Ram4]{Ra3}
D.~Ramakrishnan.
\newblock {Modularity of solvable {A}rtin representations of {${\rm
  GO}(4)$}-type}.
\newblock {\em Int. Math. Res. Not.} (2002), 1--54.

\bibitem[Ram5]{Ra4}
D.~Ramakrishnan.
\newblock {Algebraic cycles on {H}ilbert modular fourfolds and poles of
  {$L$}-functions}.
\newblock In {\em Algebraic groups and arithmetic}, pages 221--274. Tata Inst.
  Fund. Res., Mumbai, 2004.

\bibitem[Ram6]{Ra5}
D.~Ramakrishnan.
\newblock {Irreducibility of $\ell$-adic associated to regular cusp forms on
  GL$(4)/\Q$, {\rm preprint}}, 2004.

\bibitem[Ram7]{Ra6}
D.~Ramakrishnan.
\newblock {Irreducibility and Cuspidality}.
\newblock In {\em Representation Theory and Automorphic Functions}.

\bibitem[RS]{RaSh}
D.~Ramakrishnan and F.~Shahidi.
\newblock {Siegel modular forms of genus $2$ attached to elliptic curves}.
\newblock {\em Math. Res. Lett.} {\bf 14} (2007), 315--332.

\bibitem[RW]{RaW}
D.~Ramakrishnan and S.~Wang.
\newblock {A cuspidality criterion for the functorial product on {$\rm
  GL(2)\times GL(3)$} with a cohomological application}.
\newblock {\em Int. Math. Res. Not.} (2004), 1355--1394.

\bibitem[Rib]{Ri}
K.~A. Ribet.
\newblock {On {$l$}-adic representations attached to modular forms}.
\newblock {\em Invent. Math.} {\bf 28} (1975), 245--275.

\bibitem[Ser]{Se}
J.-P. Serre.
\newblock {\em Abelian {$l$}-adic representations and elliptic curves}.
\newblock Advanced Book Classics. Addison-Wesley Publishing Company Advanced
  Book Program, Redwood City, CA, second edition, 1989.
\newblock With the collaboration of Willem Kuyk and John Labute.

\bibitem[Sha1]{Sh3}
F.~Shahidi.
\newblock {On certain {$L$}-functions}.
\newblock {\em Amer. J. Math.} {\bf 103} (1981), 297--355.

\bibitem[Sha2]{Sh1}
F.~Shahidi.
\newblock {On the {R}amanujan conjecture and finiteness of poles for certain
  {$L$}-functions}.
\newblock {\em Ann. of Math. (2)} {\bf 127} (1988), 547--584.

\bibitem[Tay1]{Ta1}
R.~Taylor.
\newblock {On {G}alois representations associated to {H}ilbert modular forms}.
\newblock {\em Invent. Math.} {\bf 98} (1989), 265--280.

\bibitem[Tay2]{Ta2}
R.~Taylor.
\newblock {{$l$}-adic representations associated to modular forms over
  imaginary quadratic fields. {II}}.
\newblock {\em Invent. Math.} {\bf 116} (1994), 619--643.

\bibitem[Wan]{Wang}
S.~Wang.
\newblock {On the symmetric powers of cusp forms on {${\rm GL}(2)$} of
  icosahedral type}.
\newblock {\em Int. Math. Res. Not.} (2003), 2373--2390.

\bibitem[Wei]{W}
R.~Weissauer.
\newblock {Four dimensional {G}alois representations}.
\newblock {\em Ast\'erisque} (2005), 67--150.
\newblock Formes automorphes. II. Le cas du groupe $\rm GSp(4)$.

\end{thebibliography}

Dinakar Ramakrishnan

253-37 Caltech

Pasadena, CA 91125, USA.

dinakar@caltech.edu

\bigskip

\end{document}